\RequirePackage{fix-cm}
\documentclass[smallextended]{svjour3}       
\smartqed  
\usepackage{graphicx}
\usepackage{amsmath}
\usepackage{amssymb}
\usepackage{amsfonts}
\usepackage[ruled,vlined]{algorithm2e}
\input{tcilatex}
\begin{document}

\title{A Method to construct all the Paving Matroids over a Finite Set
\thanks{This research was partially supported by DGAPA-sabbatical-fellowship of UNAM. And by Papiit-project IN115414 of UNAM.}
}


\author{B. Mederos  \and I. P\'{e}rez-Cabrera  \and
	M. Takane \and G. Tapia-S\'{a}nchez  \and B. Zavala
}


\institute{B. Mederos and G. Tapia-S\'{a}nchez\at
              Instituto de Ingenier\'{\i}a y Tecnolog\'{\i}a \\
              Universidad Aut\'{o}noma de Ciudad Ju\'{a}rez \\
              \email{boris.mederos@uacj.mx}           
           \and
				I. P\'{e}rez-Cabrera \at
				Illumina Educaci\'{o}n \\
				\email{ivan@illumina.me}           
           \and
				M. Takane \at
				Instituto de Matem\'{a}ticas \\
				Universidad Nacional Aut\'{o}noma de M\'{e}xico (UNAM) \\
				\email{takane@matem.unam.mx (email contact)}           
			\and
				B. Zavala \at
				Facultad de Ciencias \\
				Universidad Aut\'{o}noma del Estado de M\'{e}xico \\
				\email{bzs@unam.mx}           
}

\date{Received: date / Accepted: date}

\maketitle

\begin{abstract}
We give a characterization of a matroid to be paving, through its set of hyperplanes and give an algorithm to construct all of them.
\end{abstract}


\section{Introduction}

The study of sparse-paving and paving matroids helps to understand the
behavior of matroids in general. Important examples of matroids, like the
combinatorial finite geometries are indeed paving matroids, They also have
an important role in Computer Science with the greedy algorithms and the
matroid oracles among others.

In 1959, Hartmanis \cite{H59} introduced the definition of paving matroid
through the concept of $d$-partition in number theory, see also the works of
Welsh \cite[(1976)]{We76}, Oxley \cite{O2011} and Jerrum \cite{J2006}. In
this paper, we will work with the lattice of all the subsets of a set, not
only with the so called lattice of a matroid. For references of theory of
lattices and theory of lattices of matroids, see \cite{Bi67}, \cite{H59}.

\bigskip

In this work, we give another characterization of a matroid to be paving,
which leads us to an algorithm to construct all the paving matroids. Namely,

\bigskip

\noindent \textbf{Theorem}. $\ $\textit{Let }$\ M=(S,\mathcal{I})$ be a
simple matroid of rank $r\geq 2$. Let $\mathcal{H}$ be the set of
hyperplanes of $M.$ Then

\begin{center}
	$M$ is a paving matroid if and only if $\forall X\neq Y\in \mathcal{H}$ such
	that $\left\vert X\right\vert ,\left\vert Y\right\vert \geq r.$ Then $%
	\left\vert X\cap Y\right\vert \leq r-2.$
\end{center}

\bigskip

\ The above result is a consequence of

\bigskip

\noindent \textbf{Theorem}. $\ $\textit{Let} $S=\{1,...,n\}$ \textit{be a set%
}, $r\in \mathbb{N}$ \textit{and} $2\leq r\leq n.$

\textit{Take} $\emptyset \neq \mathcal{H}^{\prime }\mathcal{\subset P}S$ 
\textit{such that} $\forall X\in \mathcal{H}^{\prime },r\leq \left\vert
X\right\vert <n$ \textit{and } $\forall X,Y\in \mathcal{H}^{\prime }$. $%
X\neq Y$ $\Longrightarrow \left\vert X\cap Y\right\vert \leq r-2.$

\textit{Let define by} $\mathcal{C}_{r}:=\{A\in \binom{S}{r};\exists X\in 
\mathcal{H}^{\prime }$ \textit{with} $A\subseteq X\}$ \textit{and} $\mathcal{%
	H}$:=$\mathcal{H}^{\prime }\cup \left\{ A\in \binom{S}{r-1};\forall D\in 
\binom{S}{r};A\subset D\text{ }\Longrightarrow \text{ }D\in \binom{S}{r}%
\backslash \mathcal{C}_{r}\right\} .$

\textit{Then} $\mathcal{H}$ \textit{\ is the set of hyperplanes of a paving
	matroid on} $S$ of rank $r$. With $\mathcal{C}_{r}$ its set of $r-$circuits
and $\mathcal{B=}\binom{S}{r}\backslash $ $\mathcal{C}_{r}$ its set of basis.

\bigskip

\noindent \textbf{Theorem}. $\ $\textit{Let} $\mathcal{M}_{n,r}$ \textit{be
	the set of the matroids on a set }$S$ \textit{with} $\left\vert S\right\vert
=n$ \textit{and rank} $r$ and $Sp_{n,t}$ be the set of sparse-paving
matroids on $S$ of rank $t.$ \textit{Then}

\begin{center}
	$\mathcal{M}_{n,r}\hookrightarrow \dprod\limits_{t=r}^{n-1}Sp_{n,t}$
\end{center}

Therefore,

\bigskip

\begin{center}
	$\left\vert \mathcal{M}_{n,r}\right\vert \leq
	\dprod\limits_{t=r}^{n-1}\left\vert Sp_{n,t}\right\vert \leq \left\vert
	Sp_{n,\frac{n}{2}}\right\vert ^{n-r}$
\end{center}

\bigskip

The material is organized as follows: In Section I, we give a
characterization of paving and sparse-paving matroids by their sets of
circuits. In section II, we give our main result: a construction of all the
paving matroids using the so called $d-$partitions. In section III, we give
an algorithm to construct the paving matroids on $S$, $\left\vert
S\right\vert =n$ and rank $r.$

\bigskip

\section{Definitions and known results}

We recall that a \textbf{matroid} $M=(S,\mathcal{I})$ consists of a finite
set $S$ and a collection $\mathcal{I}$ of subsets of $S$ (called the \textbf{%
	independent sets} of $M$) satisfying the following \textbf{independence
	axioms}:

\textit{(}$\mathcal{I}$\textit{1) }The empty set $\emptyset \in \mathcal{I}$.

\textit{(}$\mathcal{I}$\textit{2)} If $X\in \mathcal{I}$ and $Y\subseteq X$
then $Y\in \mathcal{I}$.

\textit{(}$\mathcal{I}$\textit{3)} Let $U,V\in \mathcal{I}$ \ with $%
\left\vert U\right\vert =\left\vert V\right\vert +1$ then $\exists x\in
U\backslash V$ such that $V\cup \{x\}\in \mathcal{I}$.

\medskip

A subset of $S$ which does not belong to $\mathcal{I}$ is called a \textit{%
	dependent set} of $M$.

A \textbf{basis} [respectively, a \textbf{circuit}] of $M$ is a maximal
independent [resp. minimal dependent] set of $M$.\textbf{\ }The \textbf{rank}
of a subset $X\subseteq S$ is rk$X:=\max \{\left\vert A\right\vert
;A\subseteq X$ and $A\in \mathcal{I}\}$ and the rank of the matroid $M$ is rk%
$M:=$rk$S=:r$. A \textbf{hyperplane} $X\subseteq S$ is a maximal subset of
rank $r-1.$

\noindent Known: Any matroid $M=(S,\mathcal{I})$ is completely determined by
its set of basis, $\mathcal{B}$. Namely, $\mathcal{I}=\{X\subseteq S;\exists
B\in \mathcal{B}$ with $X\subseteq B\}.$ And if $r=$rk$M$ is the rank of $M,$
then any circuit $X$ of $M,$ has cardinality $\left\vert X\right\vert \leq $%
rk$M+1$. And any hyperplane $Y$\ has cardinality $r-1\leq \left\vert
Y\right\vert \leq n-1.$

\bigskip

Let $M=(S,\mathcal{I})$ be a matroid of rank $r$. Denote by $M^{\ast }=(S,%
\mathcal{I}^{\ast })$ the \textbf{dual matroid} of $M$\textit{\ }whose set
of basis is $\mathcal{B}^{\ast }:=S\backslash \mathcal{B}$, then the rank of 
$M^{\ast }$ is $\left\vert S\right\vert -r.$

A matroid is \textbf{paving} if it has no circuits of cardinality less than
rk$M.$ And a matroid $M$ is \textbf{sparse-paving} if\ $M$ and its dual $%
M^{\ast }$ are paving matroids.

\bigskip

Along the paper, a matroid means a\textbf{\ simple matroid,} that is, it has
not circuits of cardinality 1. For general references of Theory of Matroids,
see \cite{W35}, \cite{NK2009}, \cite{We76} and \cite{O2011}.

\bigskip

\section{A description of the Paving and Sparse-paving matroids through their set of circuits}

{\LARGE \bigskip }

For any set $X$ and $m\in \mathbb{N}$, let define by $\binom{X}{m}%
=\{A\subseteq X;\left\vert A\right\vert =m\}$ the $m-$\textbf{subsets} of $X$%
.

Recall that $U_{n,r}$ is a \textbf{uniform matroid} on $S$ of cardinality $n$
and rank $r$ if $\mathcal{B}$=$\binom{S}{r}.$ That is, $\mathcal{I}%
=\{X\subseteq S;\left\vert X\right\vert \leq r\}.$ For example, any matroid
of rank $0$ or $n$ are uniform. Any uniform matroid is an sparse-paving
matroid, therefore a paving one.

For any $n$ and $r=1,$ since we work with simple matroids, it follows that
it must be $U_{n,1}$ the uniform matroid of rank $1.$

Then we will work with $2\leq r\leq n.$

\bigskip

For this section, see [10].

\noindent \textbf{Definition}. Let $M=(S,\mathcal{I})$ be a paving matroid
on $S=\{1,2,...,n\}$ of rank $2\leq r\leq n.$

Let $\mathcal{B}$ be the set of basis of $M$.

Let $\mathcal{C}_{r}$ (resp. $\mathcal{C}_{r+1}$) be the $r-$circuits ($r+1-$%
circuits), the set of the circuits of cardinality $r$ ($r+1$).

$\mathcal{N}_{1}=\{X\in \binom{S}{r+1};\exists C\in \mathcal{C}_{r}$ and $%
\exists B\in \mathcal{B}$ tal que $X=C\cup B\}.$ Moreover, for $X\in 
\mathcal{N}_{1},\exists !C\in \mathcal{C}_{r}$ such that $C\subset X.$

$\mathcal{N}_{2}=\{X\in \binom{S}{r+1};\forall A\in \binom{S}{r},A\subset X$
entonces $A\in \mathcal{C}_{r}\}.$ And $\forall X\in \mathcal{N}_{2}$, rk$%
X=r-1.$

\bigskip

\noindent \textbf{Observation}

In the lattice of \textbf{all} subsets
of $S,$ $\mathcal{L}_{S}$, we have that $\binom{S}{r-1}\subset \mathcal{I}$, 
$\binom{S}{r}=\mathcal{B}\cup \mathcal{C}_{r}$ and $\binom{S}{r+1}=\mathcal{C%
}_{r+1}\cup \mathcal{N}_{1}\cup \mathcal{N}_{2}.$

\bigskip

\noindent \textbf{Proposition[10]}\textit{\ If} $M=(S,\mathcal{I})$ 
\textit{is a paving matroid, then} $M$ \textit{is a sparse-paving if and
	only if }$\mathcal{N}_{2}=\emptyset ._{\blacksquare }$

\bigskip

Moreover, we get a method and an algorithm to construct all the
sparse-paving matroids. For possible interest, since the proof is using only
circuits, we put the proof of the next theorem in the Appendix.

\bigskip

\noindent \textbf{Theorem[10].} $\ $\textit{Let} $S=\{1,...,n\}$ 
\textit{be a set,} $r\in \mathbb{N}$ \textit{and} $r\leq n.$

\textit{Let take} $\mathcal{C}_{r}\subset \binom{S}{r}$ \textit{such that} $%
\forall C,C^{\prime }\in \binom{S}{r}$ \textit{with} $C\neq C^{\prime }$ 
\textit{then} $\left\vert C\cap C^{\prime }\right\vert \leq r-2$ \textit{and
	let} $\mathcal{B}$:=$\binom{S}{r}\backslash \mathcal{C}_{r}$ \textit{be the
	set of basis of a matroid} $M$ \textit{on} $S.$ \textit{Then} $M$ \textit{is
	a sparse-paving matroid of rank} $r._{\blacksquare }$

\bigskip

\bigskip \bigskip

\section{A description of the Paving Matroids through their set of hyperplanes}

\bigskip

\noindent \textbf{II.1. }Welsh, D. J. A. in \cite{We76} characterizes the
paving matroids in the following way:

If a paving matroid $M=(S,\mathcal{I})$ has rank $3\leq d+1<\left\vert
S\right\vert $, then its hyperplanes form a set system known as a $d$%
-partition. A family of two or more sets a $d-$\textbf{partition} if every
set in $\mathcal{F}$ has size at least $d$ and every $d$-element subset of $%
\dbigcup \mathcal{F}$ is a subset of exactly one set in $\mathcal{F}$.
Conversely, if $\mathcal{F}$ is a $d$-partition, then it can be used to
define a paving matroid on $E=\dbigcup $ $\mathcal{F}$ for which $\mathcal{F}
$ is the set of hyperplanes. See also \cite{H59}.

\bigskip

\noindent \textbf{II.2. }\noindent \textbf{Proposition}. \textit{Let }$M=(S,%
\mathcal{I})$ \textit{be a paving matroid of rank} $r\leq n$ \textit{and let 
}$\mathcal{H}$ \textit{be its set of hyperplanes.}

\textit{Then }$\mathcal{H}$\textit{\ has the following properties:}

\noindent \textbf{a)} $\forall X,Y\in \mathcal{H}$ \textit{such that} $X\neq
Y$ \textit{and} $\left\vert X\right\vert ,\left\vert Y\right\vert \geq r,$ 
\textit{we have} $\left\vert X\cap Y\right\vert \leq r-2.$

\noindent \textbf{b) }$M$ \textit{is sparse-paving if and only if} $\mathcal{%
	H\subset }\binom{S}{r-1}\cup \binom{S}{r}.$

\bigskip

\noindent \textbf{Proof:} By Welsh (II.1), $\mathcal{H}$ is an $(r-1)-$%
partition. Then $\forall A\in \binom{S}{r-1},\exists !X\in \mathcal{H}$ such
that $A\subseteq X.$ Then

\noindent \textbf{a)} Let $X,Y\in \mathcal{H}$ satisfy $X\neq Y$ and $%
\left\vert X\right\vert ,\left\vert Y\right\vert \geq r.$ Then $\binom{X}{r-1%
}\cap \binom{Y}{r-1}=\emptyset .$ Therefore, $\left\vert X\cap Y\right\vert
\leq r-2$, (otherwise, if $\left\vert X\cap Y\right\vert \geq r-1$ then $%
\exists A\subset X\cap Y$ and $\left\vert A\right\vert =r-2,$ a
contradiction).

\noindent \textbf{b) }By (I.2), $M$ is sparse-paving

if and only if $\mathcal{N}_{2}=\emptyset $

if and only if $\binom{S}{r+1}=\mathcal{C}_{r+1}\mathcal{\cup N}_{1}$ (That
is (by (I.1)), $\forall X\in \binom{S}{r+1},$ rk$X=r$)

if and only if $\mathcal{H\subset }\binom{S}{r-1}\cup \binom{S}{r}%
._{\blacksquare }$

\bigskip

\noindent \textbf{II.3. }The next result is the construction of all paving
matroids. Moreover, if there exists a hyperplane $X$ of cardinality bigger
than $r,$ then the matroid is paving no-sparse-paving.

\bigskip

\noindent \textbf{Theorem}. $\ $\textit{Let} $S=\{1,...,n\}$ \textit{be a set%
}, $r\in \mathbb{N}$ \textit{and} $2\leq r\leq n.$

\textit{Take} $\emptyset \neq \mathcal{H}^{\prime }\mathcal{\subset P}S$ 
\textit{such that} $\forall X\in \mathcal{H}^{\prime },r\leq \left\vert
X\right\vert <n$ \textit{and } $\forall X,Y\in \mathcal{H}^{\prime }$. $%
X\neq Y$ $\Longrightarrow \left\vert X\cap Y\right\vert \leq r-2.$

\textit{Let define by} $\mathcal{C}_{r}:=\{A\in \binom{S}{r};\exists X\in 
\mathcal{H}^{\prime }$ \textit{with} $A\subseteq X\}$ \textit{and} $\mathcal{%
	H}$:=$\mathcal{H}^{\prime }\cup \left\{ A\in \binom{S}{r-1};\forall D\in 
\binom{S}{r};A\subset D\text{ }\Longrightarrow \text{ }D\in \binom{S}{r}%
\backslash \mathcal{C}_{r}\right\} .$

\textit{Then} $\mathcal{H}$ \textit{\ is the set of hyperplanes of a paving
	matroid on} $S$ of rank $r$. With $\mathcal{C}_{r}$ its set of $r-$circuits
and $\mathcal{B=}\binom{S}{r}\backslash $ $\mathcal{C}_{r}$ its set of basis.

\noindent \textbf{Proof}. To prove $\mathcal{H}$ is an $(r-1)-$partition of $%
S.$ 

\noindent Define by $\mathcal{H}_{r}:=\left\{ A\in \binom{S}{r-1};\forall D\in 
\binom{S}{r};A\subset D\text{ }\Longrightarrow \text{ }D\in \binom{S}{r}%
\backslash \mathcal{C}_{r}\right\} .$

\noindent \textbf{i)} By construction $\mathcal{H}_{r-1}=\binom{S}{r-1}%
\backslash \left\{ A\in \binom{S}{r-1};\exists X\in \mathcal{H}^{\prime }%
\text{ with }A\subset X\right\} .$ Thus, $\dbigcup \mathcal{H}$=$%
\dbigcup\limits_{X\in \mathcal{H}}X=S.$

\bigskip

\noindent \textbf{ii)} To prove that $\forall A\in \binom{S}{r-1},\exists
!X\in \mathcal{H},A\subseteq X.$

\textbf{ii.a)} If $A\in \mathcal{H}_{r-1}$ we have that for all $Y\neq A$
such that $A\subset Y$, rk$Y=n.$ Therefore, $A$ is the unique hyperplane
containing $A$ itself.

\textbf{ii.b)} If $A\in \binom{S}{r-1}\backslash \mathcal{H}_{r-1}$. Then
there exists $C\in \mathcal{H}_{r}\cup \widetilde{\mathcal{C}}_{r}$ such
that $A\subset C,$ where $\widetilde{\mathcal{C}}_{r}:=\{C\in \binom{S}{r}%
;\exists X\in \mathcal{H}$, $\left\vert X\right\vert \geq r+1$ and $C\subset
X\}$ and $\mathcal{H}_{r}:=\left\{ C\in \mathcal{H};\left\vert X\right\vert
=r\right\} .$

\bigskip

\noindent Subcases: $C\in \mathcal{H}_{r}$ or $C\in \widetilde{\mathcal{C}}%
_{r}.$

\noindent \textbf{subcase(ii.b.1):} $C\in \mathcal{H}_{r}.$

By (II.2), $\forall X\in \mathcal{H}\backslash \{C\},\left\vert X\cap
C\right\vert \leq r-2,$ then (since $\left\vert A\right\vert =r-1),$ $%
A\nsubseteq X.$

\noindent \textbf{subcase(ii.b.2):} $C\in \widetilde{\mathcal{C}}_{r}.$

Then there exists $X\in \mathcal{H}$ with $\left\vert X\right\vert \geq r+1$
and $A\subset C\subset X.$ Again by (II.2), $\forall Y\in \mathcal{H}%
\backslash \{X\},$ $\left\vert A\cap Y\right\vert $ $\leq \left\vert C\cap
Y\right\vert \leq \left\vert X\cap Y\right\vert \leq r-2,$ Thus $A\nsubseteq
Y.$

\noindent Therefore,\textbf{\ }$\forall A\in \binom{S}{r-1},\exists !X\in $%
\textbf{\ }$\mathcal{H}$\textbf{\ }such that\textbf{\ }$A\subset X.$

\noindent \textbf{Therefore, }$\mathcal{H}$ is the set of hyperplanes of a
paving matroid on $S._{\blacksquare }$

\bigskip

\noindent

\noindent \textbf{Corollary. } \textit{Let }$M=(S,\mathcal{I})$ \textit{be a
	matroid of rank} $r\leq n$ \textit{and let }$\mathcal{H}$ \textit{be its set
	of hyperplanes. Then}

\begin{center}
	$M$ is a paving matroid if and only if $\forall X\neq Y\in \mathcal{H}$ such
	that $\left\vert X\right\vert ,\left\vert Y\right\vert \geq r.$ Then $%
	\left\vert X\cap Y\right\vert \leq r-2.$
	
	\bigskip
\end{center}

\section{An algorithm to construct the paving matroids}

\bigskip

Let $n,r$ be natural numbers satisfying $r\leq n-1.$

The algorithm below construct a maximal set of hyperplanes, $\mathcal{H}_{t}$%
, of cardinality $t\in \{r,...,n-1\}$ of a matroid $M=(S,\mathcal{I})$ of
rank $r.$

The hyperplanes of cardinality $r-1,$ $\mathcal{H}_{r-1}$=$\binom{S}{r-1}%
\backslash \dbigcup\limits_{t=r}^{n-1}\mathcal{H}_{t}$

\bigskip

\begin{algorithm}[H]
	\SetAlgoLined
	\KwIn{$r,n$}
	\KwOut{$H$}
	$H=\emptyset$, $k=0$\;
	$S=\{1,2, \ldots n\}$\;
	
	\While{$k<Bound$}
	{
		Choose $m \in [r,n-1]$\;
		Choose $S_i \subset S$ such that $\left\vert S_i\right\vert =m$\;
		$flag = 0$\;
		\For{all $S_j\in H$}
		{
			\If{$\left\vert S_i\cap S_j\right\vert > r-2$}{
				$flag = 1$\;
				Break\;
			}
		}
		
		\If{$flag = 0$}{
			$H=H\cup S_i$\;
		}	
		
		$k=k+1$\;
	}
	\caption{Paving Matroids}
\end{algorithm}

\bigskip

\section{Appendix: Another construction of the Sparse-paving matroids, see [10]}

\bigskip

Let $S=\{1,...,n\}$ be a set and $2\leq r\leq n.$

Let take $\mathcal{C}_{r}\subset \binom{S}{r}$ such that $\forall
C,C^{\prime }\in \mathcal{C}_{r},$ $\left\vert C\cap C^{\prime }\right\vert
\leq r-2.$ Then

\bigskip

\noindent \textbf{Theorem [10].. }\textit{Let }$S$\textit{\ be a set of
	cardinality }$\left\vert S\right\vert =n\geq 3$\textit{\ and }$2\leq r\leq
n-1.$\textit{\ Let }$C\subset \binom{S}{r}$\textit{\ be a set of }$r-$%
\textit{subsets of }$S$\textit{, satisfying the following property}

\begin{center}
	$\forall X,Y\in C$\textit{\ with }$X\neq Y$\textit{\ then }$\left\vert X\cap
	Y\right\vert \leq r-2\ \ \ \ \ \ \ \ \ (\ast \ast )$
\end{center}

\textit{Define }$M:=(S,\mathcal{I})$\textit{\ where }$\mathcal{B}:=\binom{S}{%
	r}\backslash \mathcal{C}$\textit{\ and }$\mathcal{I}:=\{X\subseteq S;$%
\textit{\ }$\exists B\in \mathcal{B}$\textit{\ with }$X\subseteq B\}.$%
\textit{\ Then, \ }\textbf{(A)}\textit{. }$M$\textit{\ is a matroid of }rk$%
M= $\textit{\ }$r$\textit{\ and }\textbf{(B)}\textit{. }$M$\textit{\ is
	sparse-paving.}

\bigskip

\noindent \textbf{Proof.} Let $S$ be a set and take a subset $\mathcal{C}%
\subset \binom{S}{r}$ satisfying the property $(\ast \ast )$. Take $M=(S,%
\mathcal{I})$ with set of basis $\mathcal{B}=\tbinom{S}{r}\backslash 
\mathcal{C}$.

\noindent \textbf{A.} \textit{To prove }$M$\textit{\ is a matroid of rank }rk%
$M=r$\textit{. }

For this proof, we will use an equivalent definition of matroid, which says:

Let $M=(S,\mathcal{I})$ is a matroid if and only if $\mathcal{I}$ satisfies (%
$\mathcal{I}$1),($\mathcal{I}$2) as in the introduction and ($\mathcal{I}$3)$%
^{\prime }$: let $B_{1},B_{2}\in \mathcal{B}$ be two basis of $M$ and $x\in
B_{1}\backslash B_{2}.$ To prove $\exists y\in B_{2}\backslash B_{1}$ such
that $(B_{1}\backslash \{x\})\cup \{y\}\in \mathcal{B}$.

\bigskip

\noindent \textbf{case a.} If $\left\vert S\right\vert =3$ and rk$M=2,$ the
possibilities for $\mathcal{C}$\ to have property $(\ast \ast )$ are $%
\mathcal{C}=\emptyset $ or $\left\vert \mathcal{C}\right\vert =1.$ In both
cases, $M$ is a matroid and it is sparse-paving.

\noindent \textbf{case b.} $\left\vert S\right\vert \geq 4.$

\textbf{(}$\mathcal{I}$\textbf{1)}\textit{\ }To prove that $\emptyset $ is
an independent set. It is enough to prove that $\mathcal{B}$ is not empty.

Since $n\geq 4$, $2\leq r\leq n-1$ and $S=\{1,...,r,r+1,...,n\}.$ Take $%
A_{1}=\{1,...,r-1,r\},$ $A_{2}=\{1,...,r-1,r+1\}$ which are subsets of $S$
with cardinality $r$ and $\left\vert A_{1}\cap A_{2}\right\vert =r-1$. Then
by $(\ast \ast ),$ $\exists i\in \{1,2\}$ such that $A_{i}\in \mathcal{B}$.
Then $\mathcal{B}\neq \emptyset .$

\textbf{(}$\mathcal{I}$\textbf{2)} Let $Y\subseteq X\subseteq S$ such that $%
\exists B\in \mathcal{B}$ with $X\subseteq B$. Then $Y\subseteq B$, that is $%
Y$ is independent, by definition.

\textbf{(}$\mathcal{I}$\textbf{3)}$^{\prime }$ Now, let $B_{1},B_{2}\in 
\mathcal{B}$ be two basis of $M$ and $x\in B_{1}\backslash B_{2}.$ To prove $%
\exists y\in B_{2}\backslash B_{1}$ such that $(B_{1}\backslash \{x\})\cup
\{y\}\in \mathcal{B}$.

$\mathcal{I}$\textbf{3}$^{\prime }$\textbf{.1.} Assume $m:=\left\vert
B_{2}\backslash B_{1}\right\vert =1$. That is, $B_{2}\cap
B_{1}=B_{1}\backslash \{x\}$ and $B_{2}=(B_{1}\backslash \{x\})\cup \{y\}$
for some $y\in S.$ Then $(B_{1}\backslash \{x\})\cup \{y\}\in \mathcal{B}$.

$\mathcal{I}$\textbf{3}$^{\prime }$\textbf{.2.} Let define $m:=\left\vert
B_{2}\backslash B_{1}\right\vert \geq 2$ and let $B_{2}=$ $(B_{1}\cap
B_{2})\cup \{y_{1},y_{2},y_{3},...,y_{m}\}$. Define $A_{i}:=(B_{1}\backslash
\{x\})\cup \{y_{i}\}$ for $i=1,...,m.$ Since $\forall i\neq j,$ $\left\vert
A_{i}\cap A_{j}\right\vert =r-1$ and $m\geq 2$, by $(\ast \ast ),$ $\exists
A_{i_{0}}\in \mathcal{B}$. Therefore, $(B_{1}\backslash \{x\})\cup
\{y_{i}\}=A_{i_{0}}\in \mathcal{B},$ and $M$ \textit{is a matroid}.

\emph{Rank: }By definition of $M$, rk$M=r$.

\bigskip

\noindent \textbf{B.}\textit{\ To prove }$M$\textit{\ is a sparse-paving
	matroid. }

\textbf{B.1.} First we will prove that $M$ is a \textit{paving matroid}.
Equivalently, to prove $\forall Z\subseteq S$ of $\left\vert Z\right\vert =$%
rk$M-1,$ $Z\in \mathcal{I}$. $\ $This proof is similar to the one of ($%
\mathcal{I}$1). Namely:

Let rk$M\leq n-1$. Since $n\geq 3$ and $\left\vert Z\right\vert =$rk$M-1$,
we have $S=Z\cup \{x_{1},x_{2},...,x_{m}\}$ with $m\geq 2$. Let denote $%
A_{i}:=Z\cup \{x_{i}\}$ for $i=1,2,...,m.$ By $(\ast \ast )$ and $m\geq 2,$ $%
\exists i_{0}\in \{1,...,m\}$ such that $(Z\subset )A_{i_{0}}\in \mathcal{B}$%
. Then $Z\in \mathcal{I}$.

\textbf{B.2.} And by (1.2), $M$ is a sparse-paving matroid.$_{\blacksquare }$

\bigskip



%
%



\end{document}